\documentclass[12pt,a4paper]{article}

\usepackage{amsmath,amssymb,amsthm}
\usepackage{amsfonts}
\usepackage[utf8]{inputenc}
\usepackage[margin=2.5cm]{geometry}
\newtheorem{theorem}{Theorem}[section]
\newtheorem{lemma}[theorem]{Lemma}
\newtheorem{proposition}[theorem]{Proposition}

\theoremstyle{definition}
\newtheorem{definition}[theorem]{Definition}
\newtheorem{example}[theorem]{Example}

\numberwithin{equation}{section}

\begin{document}

\title{\textbf{Soft $g$-Frames in Soft Hilbert Spaces}}
\author{\textbf{Sayyed Mehrab Ramezani}}
\date{}

\maketitle

\begin{abstract}
In this paper, we introduce and analyze soft $g$-frames in soft Hilbert spaces, extending generalized frame theory to parameterized uncertainty models. We establish the completeness of the direct sum of soft Hilbert spaces and prove that the associated soft pre-$g$-frame operator $T_\Lambda$ is well-defined and bounded with unconditionally convergent series representations. Furthermore, we show that the soft $g$-frame operator $S_\Lambda$ is a bounded, self-adjoint, and strictly positive element of the unital soft Banach algebra $\mathcal{B}(\tilde{\mathcal{U}})$,
and we prove its invertibility via the soft Neumann series and derive two-sided operator bounds for its inverse. We characterize the canonical dual soft $g$-frame, establish the unconditional soft $g$-frame decomposition theorem, and construct induced pairs of dual soft frames from dual soft $g$-frames. Finally, an analytical example with parameter-dependent metrics illustrates the matrix representations, verifies the optimality of the frame bounds via Rayleigh quotients, and confirms the reconstruction formula.
\end{abstract}

\vspace{1ex}

\vspace{1ex}

\section{Introduction}
Frame theory in Hilbert spaces was introduced by Duffin and Schaeffer \cite{duf} in the context of non-harmonic Fourier series. Frames provide redundant, stable representations of vectors, overcoming the strict linear independence constraints of orthonormal bases. To unify various generalized frame concepts---including fusion frames, pseudo-frames, and bounded linear operator representations---Sun \cite{sun2006, sun} formulated the theory of $g$-frames and $g$-Riesz bases in Hilbert spaces, establishing their fundamental operator properties, perturbation bounds, and duality characterizations. Expanding this framework, Ramezani and Nazari \cite{meh} examined structural generalizations including $g$-orthonormal bases, $g$-Riesz bases, and the characterization of $g$-duals.

In parallel with developments in functional analysis, Molodtsov \cite{mol} initiated soft set theory as a parameterized mathematical tool to handle uncertainty, vagueness, and incomplete data without the sharp boundary limitations of classical models. Maji et al. \cite{maj, roy} analyzed theoretical operations on soft sets and demonstrated their applicability to decision-making problems. Subsequently, algebraic and analytic structures were established on soft sets. Das and Samanta introduced soft real sets and soft real numbers \cite{das1}, soft complex sets and numbers \cite{das2}, soft metric spaces \cite{das}, and soft normed linear spaces with associated bounded soft linear operators \cite{das4}. Furthermore, Das and Samanta developed the foundation of soft inner product spaces and established operator theory on soft Hilbert spaces \cite{das5, das3}. Structural completeness and spectral properties were later generalized to soft Banach algebras by Thakur and Samanta \cite{Tak}.

Integrating redundant representations with soft structures, Ferrer, Sierra, and Sanabria \cite{soft} introduced discrete soft frames in soft Hilbert spaces. They established the boundedness of the soft pre-frame operator, demonstrated the invertibility of the soft frame operator using the Neumann series in soft Banach algebras, and proved the soft frame decomposition theorem.

\subsection{Structural Perspective and Relation to Classical $g$-Frames}

While the formulation in \cite{soft} addresses discrete vector systems in soft Hilbert spaces, generalized operator-based representations within soft structures remain undeveloped. To clarify the mathematical scope of the present framework, it is important to examine the relationship between soft $g$-frames and families of classical $g$-frames \cite{sun}.

When the parameter set $A$ is finite, the algebraic operations, soft norms, and partial order relations in $\mathbb{R}(A)$ are defined pointwise across $A$. Consequently, by the decomposition properties of soft inner product spaces \cite{das5, das3}, a soft $g$-frame $\{\Lambda_j\}_{j\in J}$ for $\tilde{\mathcal{U}}$ with respect to $\{\tilde{\mathcal{V}}_j\}_{j\in J}$ decomposes pointwise into a parameterized family of classical $g$-frames $\{\Lambda_{j,\lambda}\}_{j\in J}$ for $\mathcal{U}_\lambda$ with respect to $\{\mathcal{V}_{j,\lambda}\}_{j\in J}$ for each $\lambda \in A$, governed by parameter-dependent bounds $\overline{C}(\lambda)$ and $\overline{D}(\lambda)$.

To avoid treating each parameter regime in isolation, formulating this theory directly within the soft Hilbert space and soft Banach algebra framework provides a unified operator-theoretic formulation across the entire parameter domain:
\begin{enumerate}
    \item \textbf{Unified Operator-Algebraic Formulation:} Rather than performing component-wise operator inversions separately for each parameter state, the soft $g$-frame operator $S_\Lambda$ and its canonical dual $\{\widetilde{\Lambda}_j\}_{j\in J}$ are formulated and inverted within a single algebraic expression in the unital soft Banach algebra $\mathcal{B}(\tilde{\mathcal{U}})$ via the soft Neumann series.
    \item \textbf{Parameter-Dependent Frame Bounds:} Although the finiteness of $A$ guarantees the existence of uniform scalar bounds $c = \min_{\lambda\in A}\overline{C}(\lambda) > 0$ and $d = \max_{\lambda\in A}\overline{D}(\lambda) < \infty$, retaining soft real bounds $\overline{C}, \overline{D} \in \mathbb{R}(A)^*_{>0}$ preserves sharp, parameter-specific stability margins across different environmental regimes $\lambda \in A$ within a unified reconstruction identity.
    \item \textbf{Outlook toward Infinite Parameter Spaces:} Establishing the direct sum topology, unconditional convergence criteria, and operator-algebraic duality in the finite-parameter setting provides a natural starting point for future investigations on infinite or topological parameter sets, where finite-product decompositions are no longer available.
\end{enumerate}
In summary, the primary contribution of this work is not to assert a stronger analytical theory than classical $g$-frame families in the finite-parameter case, but rather to establish a cohesive, unified operator-theoretic calculus in the language of soft Hilbert spaces and soft Banach algebras.
\section{Preliminaries and Definitions}

In this section, we review foundational concepts concerning soft sets, soft Hilbert spaces, and bounded soft linear operators. Throughout this paper, let $X$ denote a non-empty set, $\mathcal{P}(X)$ the power set of $X$, and $A$ a non-empty finite set of parameters.

\subsection{Soft Sets and Soft Real Numbers}

\begin{definition}[\cite{mol}]\label{def:soft-set}
A pair $(F, A)$ is called a \textbf{soft set} over $X$ if $F$ is a mapping given by
\[
F: A \longrightarrow \mathcal{P}(X).
\]
Equivalently, $(F, A)$ can be represented by its parameterized graph $\{( \lambda, F(\lambda) ) : \lambda \in A, \, F(\lambda) \in \mathcal{P}(X)\}$.
\end{definition}

A soft set $(F, A)$ over $X$ is termed a \emph{null soft set}, denoted by $\Phi$, if $F(\lambda) = \emptyset$ for all $\lambda \in A$. It is termed an \emph{absolute soft set}, denoted by $\tilde{X}$, if $F(\lambda) = X$ for all $\lambda \in A$.

\begin{definition}[\cite{das}]\label{def:soft-element}
A mapping $\epsilon : A \to X$ is called a \textbf{soft element} of $X$. If $\epsilon(\lambda) \in F(\lambda)$ for all $\lambda \in A$, then $\epsilon$ is said to belong to the soft set $(F, A)$, denoted by $\epsilon \tilde{\in} (F, A)$. The collection of all soft elements belonging to $(F, A)$ is denoted by
\[
SE((F, A)) := \{ \epsilon : A \to X \mid \epsilon(\lambda) \in F(\lambda), \, \forall \lambda \in A \}.
\]
\end{definition}

\begin{definition}[\cite{das1, das2}]\label{def:soft-numbers}
Let $\mathbb{K}$ denote either $\mathbb{R}$ or $\mathbb{C}$. A mapping $F: A \to \mathcal{P}(\mathbb{K})$ such that $F(\lambda)$ is a non-empty bounded subset of $\mathbb{K}$ for each $\lambda \in A$ is called a \emph{soft $\mathbb{K}$-set}. If $F(\lambda)$ is a singleton for every $\lambda \in A$, it is identified with a mapping $\alpha: A \to \mathbb{K}$ and termed a \textbf{soft $\mathbb{K}$-number}.

The set of all soft real numbers is denoted by $\mathbb{R}(A)$, and the set of all soft complex numbers is denoted by $\mathbb{C}(A)$. A constant soft number $\alpha(\lambda) = c$ for all $\lambda \in A$ is denoted by $\overline{c}$.
\end{definition}

\begin{definition}[\cite{das1}]\label{def:soft-order}
The collection of non-negative soft real numbers is denoted by
\[
\mathbb{R}(A)^* := \{ \alpha \in \mathbb{R}(A) : \alpha(\lambda) \ge 0, \, \forall \lambda \in A \}.
\]
Furthermore, the collection of strictly positive soft real numbers is denoted by
\[
\mathbb{R}(A)^*_{>0} := \{ \alpha \in \mathbb{R}(A) : \alpha(\lambda) > 0, \, \forall \lambda \in A \}.
\]
For two soft real numbers $\alpha, \beta \in \mathbb{R}(A)$, the partial order relations are defined pointwise:
\begin{enumerate}
\item $\alpha \tilde{\le} \beta$ if and only if $\alpha(\lambda) \le \beta(\lambda)$ for all $\lambda \in A$,
\item $\alpha \tilde{<} \beta$ if and only if $\alpha(\lambda) < \beta(\lambda)$ for all $\lambda \in A$.
\end{enumerate}
\end{definition}

\subsection{Soft Hilbert Spaces}

\begin{definition}[\cite{das5, das4}]\label{def:soft-vector-space}
Let $X$ be a vector space over $\mathbb{K}$. A soft set $(F, A)$ over $X$ is called a \textbf{soft vector space} if $F(\lambda)$ is a vector subspace of $X$ for each $\lambda \in A$. The algebraic operations on $SE((F, A))$ are defined pointwise across the parameter set $A$:
\[
(x + y)(\lambda) = x(\lambda) + y(\lambda), \quad (k \cdot x)(\lambda) = k(\lambda) x(\lambda), \quad \forall \lambda \in A,
\]
for all $x, y \in SE((F, A))$ and soft scalars $k: A \to \mathbb{K}$. The null soft vector $\Theta$ is defined by $\Theta(\lambda) = \theta$ for all $\lambda \in A$, where $\theta$ is the zero vector of $X$.
\end{definition}

\begin{definition}[\cite{das5}]\label{def:soft-inner-product}
Let $\tilde{X}$ be an absolute soft vector space over $\mathbb{K}$. A mapping $\langle \cdot, \cdot \rangle : SE(\tilde{X}) \times SE(\tilde{X}) \to \mathbb{C}(A)$ is called a \textbf{soft inner product} on $\tilde{X}$ if it satisfies:
\begin{enumerate}
\item $\langle x, x \rangle \tilde{\ge} \overline{0}$ for all $x \tilde{\in} \tilde{X}$, and $\langle x, x \rangle = \overline{0} \iff x = \Theta$;
\item $\langle x, y \rangle = \overline{\langle y, x \rangle}$ for all $x, y \tilde{\in} \tilde{X}$ (where the bar denotes complex conjugation);
\item $\langle \alpha \cdot x, y \rangle = \alpha \cdot \langle x, y \rangle$ for all $x, y \tilde{\in} \tilde{X}$ and soft scalars $\alpha \in \mathbb{C}(A)$;
\item $\langle x + y, z \rangle = \langle x, z \rangle + \langle y, z \rangle$ for all $x, y, z \tilde{\in} \tilde{X}$.
\end{enumerate}
The pair $(\tilde{X}, \langle \cdot, \cdot \rangle, A)$ is called a \textbf{soft inner product space}. The induced soft norm is defined by $\|x\| := \sqrt{\langle x, x \rangle} \in \mathbb{R}(A)^*$. If $\tilde{X}$ is complete with respect to the soft metric induced by this soft norm, it is called a \textbf{soft Hilbert space}.
\end{definition}

\begin{theorem}[\cite{das5, soft}]\label{thm:pointwise-inner-product}
Let $H$ be a Hilbert space and $A$ a finite parameter set. The mapping $\langle \cdot, \cdot \rangle : SE(\tilde{H}) \times SE(\tilde{H}) \to \mathbb{C}(A)$ given by
\[
\langle h_1, h_2 \rangle(\lambda) := \langle h_1(\lambda), h_2(\lambda) \rangle_{H}, \quad \forall \lambda \in A,
\]
defines a soft inner product, and $(\tilde{H}, \langle \cdot, \cdot \rangle, A)$ is a soft Hilbert space.
\end{theorem}

\subsection{Soft Linear Operators}

\begin{definition}[\cite{das3, das4}]\label{def:soft-operators}
Let $\tilde{X}$ and $\tilde{Y}$ be soft normed linear spaces. A mapping $T: SE(\tilde{X}) \to SE(\tilde{Y})$ is said to be a \textbf{soft linear operator} if:
\begin{enumerate}
\item $T(x_1 + x_2) = T(x_1) + T(x_2)$ for all $x_1, x_2 \in SE(\tilde{X})$,
\item $T(c \cdot x) = c \cdot T(x)$ for all soft scalars $c$ and $x \in SE(\tilde{X})$.
\end{enumerate}
An operator $T$ is said to be \textbf{bounded} if there exists a non-negative soft real number $\overline{M} \in \mathbb{R}(A)^*$ such that
\[
\|T(x)\| \tilde{\le} \overline{M} \|x\|, \quad \forall x \in SE(\tilde{X}).
\]
The collection of all bounded soft linear operators from $\tilde{X}$ to $\tilde{Y}$ is denoted by $\mathcal{B}(\tilde{X}, \tilde{Y})$.
\end{definition}

\begin{definition}[\cite{das3}]\label{def:soft-adjoint}
Let $\tilde{H}_1$ and $\tilde{H}_2$ be soft Hilbert spaces and $T \in \mathcal{B}(\tilde{H}_1, \tilde{H}_2)$. The unique operator $T^* \in \mathcal{B}(\tilde{H}_2, \tilde{H}_1)$ satisfying
\[
\langle T x, y \rangle_{\tilde{H}_2} = \langle x, T^* y \rangle_{\tilde{H}_1}, \quad \forall x \in SE(\tilde{H}_1), \, y \in SE(\tilde{H}_2),
\]
is called the \textbf{adjoint soft operator} of $T$. An operator $T \in \mathcal{B}(\tilde{H}, \tilde{H})$ is called \textbf{self-adjoint} if $T^* = T$, and \textbf{positive} (denoted by $T \tilde{\ge} \overline{0} Id_{\tilde{H}}$) if $\langle T x, x \rangle_{\tilde{H}} \tilde{\ge} \overline{0}$ for all $x \in SE(\tilde{H})$.
\end{definition}
\begin{theorem}[\cite{das3}]\label{thm:adjoint-existence}
Let $\tilde{H}_1$ and $\tilde{H}_2$ be soft Hilbert spaces over a finite parameter set $A$, and let $T \in \mathcal{B}(\tilde{H}_1, \tilde{H}_2)$. Then there exists a unique bounded soft linear operator $T^* \in \mathcal{B}(\tilde{H}_2, \tilde{H}_1)$ such that
\[
\langle T x, y \rangle_{\tilde{H}_2} = \langle x, T^* y \rangle_{\tilde{H}_1}, \quad \forall x \in SE(\tilde{H}_1), \, y \in SE(\tilde{H}_2).
\]
Moreover, $\|T^*\| = \|T\|$.
\end{theorem}

\begin{proposition}[\cite{das4, Tak}]\label{prop:soft-banach-algebra}
Let $\tilde{\mathcal{U}}$ be a soft Hilbert space over a finite parameter set $A$. Equipped with pointwise algebraic operations, operator composition, the identity operator $Id_{\tilde{\mathcal{U}}}$, and the soft operator norm
\[
\|T\|(\lambda) := \sup_{\|x(\lambda)\|_{\mathcal{U}_\lambda} = 1} \|(Tx)(\lambda)\|_{\mathcal{U}_\lambda}, \quad \forall \lambda \in A,
\]
the space $\mathcal{B}(\tilde{\mathcal{U}}) := \mathcal{B}(\tilde{\mathcal{U}}, \tilde{\mathcal{U}})$ forms a unital soft Banach algebra with unit $Id_{\tilde{\mathcal{U}}}$, satisfying $\|Id_{\tilde{\mathcal{U}}}\| = \overline{1}$ and the submultiplicative inequality $\|S T\| \tilde{\le} \|S\| \|T\|$ for all $S, T \in \mathcal{B}(\tilde{\mathcal{U}})$. Since $A$ is finite and each component $\mathcal{B}(\mathcal{U}_\lambda)$ is complete, $\mathcal{B}(\tilde{\mathcal{U}})$ is complete with respect to the induced soft metric.
\end{proposition}
\begin{proof}
For each $\mu \in A$, define the soft characteristic scalar $e_\mu \in \mathbb{R}(A)$ by $e_\mu(\mu) = 1$ and $e_\mu(\lambda) = 0$ for $\lambda \neq \mu$. By soft linearity of any $T \in \mathcal{B}(\tilde{\mathcal{U}})$, we have $(T(e_\mu \cdot x))(\mu) = (e_\mu \cdot Tx)(\mu) = (Tx)(\mu)$, which shows that $(Tx)(\mu)$ depends solely on $x(\mu)$. Thus, each $T$ decomposes uniquely into component operators $T_\lambda \in \mathcal{B}(\mathcal{U}_\lambda)$ via $(Tx)(\lambda) = T_\lambda x(\lambda)$.

Consequently, the canonical map $\Phi : \mathcal{B}(\tilde{\mathcal{U}}) \to \prod_{\lambda\in A} \mathcal{B}(\mathcal{U}_\lambda)$ given by $\Phi(T) := \{T_\lambda\}_{\lambda\in A}$ is a bijective soft algebra homomorphism preserving operator composition ($\Phi(ST) = \{S_\lambda T_\lambda\}_{\lambda\in A}$), mapping $Id_{\tilde{\mathcal{U}}}$ to $\{Id_{\mathcal{U}_\lambda}\}_{\lambda\in A}$, and satisfying $\|T\|(\lambda) = \|T_\lambda\|$. In particular, the submultiplicative property holds pointwise:
\[
\|ST\|(\lambda) = \|S_\lambda T_\lambda\| \le \|S_\lambda\|\|T_\lambda\| = \|S\|(\lambda)\|T\|(\lambda), \quad \forall \lambda \in A.
\]
Since $A$ is finite and each component $\mathcal{B}(\mathcal{U}_\lambda)$ is a classical Banach algebra, any Cauchy sequence $\{T^{(n)}\}_{n=1}^\infty$ in $\mathcal{B}(\tilde{\mathcal{U}})$ corresponds to Cauchy sequences $\{T^{(n)}_\lambda\}_{n=1}^\infty$ in $\mathcal{B}(\mathcal{U}_\lambda)$ converging to limits $T_\lambda \in \mathcal{B}(\mathcal{U}_\lambda)$. Defining $T := \Phi^{-1}(\{T_\lambda\}_{\lambda\in A})$ yields $\|T^{(n)} - T\| \to \overline{0}$, establishing completeness.
\end{proof}

\begin{theorem}[\cite{Tak}]\label{thm:neumann}
Let $\mathfrak{A}$ be a unital soft Banach algebra with unit $\overline{e}$. If $x \in \mathfrak{A}$ satisfies $\|\overline{e} - x\| \tilde{<} \overline{1}$, then $x$ is invertible in $\mathfrak{A}$ and its inverse is given by the convergent Neumann series:
\[
x^{-1} = \sum_{n=0}^{\infty} (\overline{e} - x)^n.
\]
\end{theorem}

\section{Soft $g$-Frames in Soft Hilbert Spaces}

Throughout this section, let $A$ be a non-empty finite set of parameters, $\tilde{\mathcal{U}} = (\mathcal{U}, A)$ be a soft Hilbert space over $\mathbb{K}$ ($\mathbb{K} = \mathbb{R}$ or $\mathbb{C}$), and $\{\tilde{\mathcal{V}}_j\}_{j\in J}$ be a family of soft Hilbert spaces parameterized by the same set $A$, where the index set $J$ is at most countable.

\subsection{Direct Sums and Convergence in Soft Hilbert Spaces}

Throughout this paper, convergence of series in $\mathbb{R}(A)$ or in a soft Hilbert space $\tilde{\mathcal{H}}$ is understood with respect to the soft metric induced by the soft norm; since $A$ is finite, this is equivalent to pointwise convergence for every $\lambda \in A$. Specifically, a countable series $\sum_{j\in J} x_j$ in $\tilde{\mathcal{H}}$ is said to \emph{converge unconditionally} to $x \in SE(\tilde{\mathcal{H}})$ if for every $\lambda \in A$, the net of finite partial sums $\{\sum_{j\in F} x_j(\lambda)\}_{F \subset J, \, |F|<\infty}$ converges to $x(\lambda)$ in $\mathcal{H}_\lambda$ independently of the ordering of $J$.

To establish a rigorous operator domain, we define the direct sum of the family $\{\tilde{\mathcal{V}}_j\}_{j\in J}$ as
\[
\bigoplus_{j\in J} \tilde{\mathcal{V}}_j := \left\{ g = \{g_j\}_{j\in J} : g_j \in SE(\tilde{\mathcal{V}}_j), \, \sum_{j\in J} \|g_j(\lambda)\|_{\mathcal{V}_{j,\lambda}}^2 < \infty \text{ for all } \lambda \in A \right\}.
\]
Equivalently, $\sum_{j\in J} \|g_j\|^2$ defines a well-defined element in $\mathbb{R}(A)^*$. For any $g = \{g_j\}_{j\in J}, h = \{h_j\}_{j\in J} \in \bigoplus_{j\in J} \tilde{\mathcal{V}}_j$, we define
\[
\langle g, h \rangle_{\oplus}(\lambda) := \sum_{j\in J} \langle g_j(\lambda), h_j(\lambda) \rangle_{\mathcal{V}_{j, \lambda}}, \quad \forall \lambda \in A,
\]
with induced soft norm $\|g\|_\oplus(\lambda) := \left( \sum_{j\in J} \|g_j(\lambda)\|_{\mathcal{V}_{j,\lambda}}^2 \right)^{1/2}$.

\begin{proposition}\label{prop:direct-sum-hilbert}
Equipped with pointwise vector space operations and the soft inner product $\langle \cdot, \cdot \rangle_\oplus$, the direct sum $\bigoplus_{j\in J} \tilde{\mathcal{V}}_j$ is a soft Hilbert space.
\end{proposition}

\begin{proof}
For each $\lambda \in A$, the component space $\bigoplus_{j\in J} \mathcal{V}_{j,\lambda}$ is the classical $\ell^2$-direct sum of the Hilbert spaces $\{\mathcal{V}_{j,\lambda}\}_{j\in J}$, which is a complete Hilbert space under the inner product $\langle \cdot, \cdot \rangle_{\oplus}(\lambda)$. Since $A$ is finite, any Cauchy sequence $\{g^{(n)}\}_{n=1}^\infty$ in $\bigoplus_{j\in J} \tilde{\mathcal{V}}_j$ satisfies $\|g^{(n)} - g^{(m)}\|_\oplus(\lambda) \to 0$ as $n,m \to \infty$ for every $\lambda \in A$. Hence, $\{g^{(n)}(\lambda)\}_{n=1}^\infty$ converges to a unique element $g(\lambda) \in \bigoplus_{j\in J} \mathcal{V}_{j,\lambda}$ for each $\lambda \in A$. Defining $g(\lambda) := \{g_j(\lambda)\}_{j\in J}$ yields an element $g \in \bigoplus_{j\in J} \tilde{\mathcal{V}}_j$ such that $\|g^{(n)} - g\|_\oplus \to \overline{0}$, proving completeness.
\end{proof}

\begin{definition}\label{def:soft-g-frame}
A sequence of bounded soft linear operators $\{\Lambda_j \in \mathcal{B}(\tilde{\mathcal{U}}, \tilde{\mathcal{V}}_j) : j \in J\}$ is called a \textbf{soft $g$-frame} for $\tilde{\mathcal{U}}$ with respect to $\{\tilde{\mathcal{V}}_j\}_{j\in J}$ if there exist soft real numbers $\overline{C}, \overline{D} \in \mathbb{R}(A)^*_{>0}$ with $\overline{C} \tilde{\le} \overline{D}$ such that
\begin{equation}\label{eq:soft-g-frame-ineq}
\overline{C} \|f\|^2 \tilde{\le} \sum_{j\in J} \|\Lambda_j f\|^2 \tilde{\le} \overline{D} \|f\|^2, \quad \forall f \in SE(\tilde{\mathcal{U}}).
\end{equation}
The elements $\overline{C}$ and $\overline{D}$ are called the lower and upper bounds of the soft $g$-frame, respectively. If $\overline{C} = \overline{D}$, the sequence is called a \textbf{tight soft $g$-frame}, and if $\overline{C} = \overline{D} = \overline{1}$, it is called a \textbf{Parseval soft $g$-frame}.
\end{definition}

\subsection{The Pre-$g$-Frame and Soft $g$-Frame Operators}

\begin{definition}\label{def:pre-g-op}
Let $\{\Lambda_j\}_{j\in J}$ be a soft $g$-frame for $\tilde{\mathcal{U}}$ with respect to $\{\tilde{\mathcal{V}}_j\}_{j\in J}$. The \textbf{soft pre-$g$-frame operator} $T_\Lambda : \bigoplus_{j\in J} \tilde{\mathcal{V}}_j \to SE(\tilde{\mathcal{U}})$ is formally defined by
\[
T_\Lambda(\{g_j\}_{j\in J}) := \sum_{j\in J} \Lambda_j^* g_j,
\]
whose unconditional convergence in the soft Hilbert space metric is established in Proposition~\ref{prop:T-bounded}.
\end{definition}

\begin{proposition}\label{prop:T-bounded}
The soft pre-$g$-frame operator $T_\Lambda$ is well-defined, soft linear, and bounded with $\|T_\Lambda\| \tilde{\le} \sqrt{\overline{D}}$. Furthermore, its adjoint soft operator $T_\Lambda^* : SE(\tilde{\mathcal{U}}) \to \bigoplus_{j\in J} \tilde{\mathcal{V}}_j$ is given by
\[
T_\Lambda^*(f) = \{\Lambda_j f\}_{j\in J}, \quad \forall f \in SE(\tilde{\mathcal{U}}).
\]
\end{proposition}

\begin{proof}
Let $g = \{g_j\}_{j\in J} \in \bigoplus_{j\in J} \tilde{\mathcal{V}}_j$ and fix an arbitrary parameter $\lambda \in A$. For any finite subset $F \subset J$ and any vector $u \in \mathcal{U}_\lambda$, applying the classical Cauchy-Schwarz inequality in $\mathcal{V}_{j,\lambda}$ and $\ell^2(F)$ yields
\begin{align*}
\left| \left\langle \sum_{j\in F} \Lambda_{j,\lambda}^* g_j(\lambda), u \right\rangle_{\mathcal{U}_\lambda} \right| 
&= \left| \sum_{j\in F} \langle g_j(\lambda), \Lambda_{j,\lambda} u \rangle_{\mathcal{V}_{j,\lambda}} \right| \\
&\le \sum_{j\in F} \|g_j(\lambda)\|_{\mathcal{V}_{j,\lambda}} \|\Lambda_{j,\lambda} u\|_{\mathcal{V}_{j,\lambda}} \\
&\le \left( \sum_{j\in F} \|g_j(\lambda)\|_{\mathcal{V}_{j,\lambda}}^2 \right)^{1/2} \left( \sum_{j\in F} \|\Lambda_{j,\lambda} u\|_{\mathcal{V}_{j,\lambda}}^2 \right)^{1/2} \\
&\le \sqrt{\overline{D}(\lambda)} \left( \sum_{j\in F} \|g_j(\lambda)\|_{\mathcal{V}_{j,\lambda}}^2 \right)^{1/2} \|u\|_{\mathcal{U}_\lambda}.
\end{align*}

Taking the supremum over all $u \in \mathcal{U}_\lambda$ with $\|u\|_{\mathcal{U}_\lambda} = 1$, we obtain for every finite subset $F \subset J$:
\[
\left\| \sum_{j\in F} \Lambda_{j,\lambda}^* g_j(\lambda) \right\|_{\mathcal{U}_\lambda} \le \sqrt{\overline{D}(\lambda)} \left( \sum_{j\in F} \|g_j(\lambda)\|_{\mathcal{V}_{j,\lambda}}^2 \right)^{1/2}.
\]
In particular, for any two finite subsets $F_1 \subset F_2 \subset J$, applying this bound to the finite difference $F_2 \setminus F_1$ yields
\[
\left\| \sum_{j\in F_2} \Lambda_{j,\lambda}^* g_j(\lambda) - \sum_{j\in F_1} \Lambda_{j,\lambda}^* g_j(\lambda) \right\|_{\mathcal{U}_\lambda} = \left\| \sum_{j\in F_2 \setminus F_1} \Lambda_{j,\lambda}^* g_j(\lambda) \right\|_{\mathcal{U}_\lambda} \le \sqrt{\overline{D}(\lambda)} \left( \sum_{j\in F_2 \setminus F_1} \|g_j(\lambda)\|_{\mathcal{V}_{j,\lambda}}^2 \right)^{1/2}.
\]
Since $\sum_{j\in J} \|g_j(\lambda)\|_{\mathcal{V}_{j,\lambda}}^2 < \infty$, the tail sum over $F_2 \setminus F_1$ tends to zero as $F_1, F_2$ expand through the directed set of finite subsets of $J$. Thus, the net of partial sums $\{\sum_{j\in F} \Lambda_{j,\lambda}^* g_j(\lambda)\}_{F \subset J, |F|<\infty}$ satisfies the Cauchy criterion in the Hilbert space $\mathcal{U}_\lambda$. Therefore, the series $\sum_{j\in J} \Lambda_j^* g_j$ converges unconditionally in $\tilde{\mathcal{U}}$, proving that $T_\Lambda$ is well-defined and soft linear.

 Passing to the limit over all finite subsets $F \subset J$ gives
\[
\|(T_\Lambda g)(\lambda)\|_{\mathcal{U}_\lambda} \le \sqrt{\overline{D}(\lambda)} \|g\|_\oplus(\lambda), \quad \forall \lambda \in A,
\]
which confirms that $T_\Lambda$ is bounded with $\|T_\Lambda\| \tilde{\le} \sqrt{\overline{D}}$.

For the adjoint operator, for any $f \in SE(\tilde{\mathcal{U}})$ and $g \in \bigoplus_{j\in J} \tilde{\mathcal{V}}_j$, continuity of the inner product and unconditional convergence on each $\lambda \in A$ give
\[
\langle T_\Lambda(g), f \rangle_{\tilde{\mathcal{U}}}(\lambda) = \sum_{j\in J} \langle \Lambda_{j,\lambda}^* g_j(\lambda), f(\lambda) \rangle_{\mathcal{U}_\lambda} = \sum_{j\in J} \langle g_j(\lambda), \Lambda_j f(\lambda) \rangle_{\mathcal{V}_{j,\lambda}} = \langle g, \{\Lambda_j f\}_{j\in J} \rangle_{\oplus}(\lambda).
\]
Hence, $T_\Lambda^*(f) = \{\Lambda_j f\}_{j\in J}$.
\end{proof}

\subsection{Invertibility and Canonical Duality}

\begin{definition}\label{def:frame-operator}
Let $\{\Lambda_j\}_{j\in J}$ be a soft $g$-frame for $\tilde{\mathcal{U}}$ with respect to $\{\tilde{\mathcal{V}}_j\}_{j\in J}$. The \textbf{soft $g$-frame operator} $S_\Lambda : SE(\tilde{\mathcal{U}}) \to SE(\tilde{\mathcal{U}})$ is defined by
\[
S_\Lambda(f) := T_\Lambda T_\Lambda^*(f) = \sum_{j\in J} \Lambda_j^* \Lambda_j f, \quad \forall f \in SE(\tilde{\mathcal{U}}).
\]
\end{definition}

\begin{definition}\label{def:duals}
Let $\{\Lambda_j\}_{j\in J}$ and $\{\Gamma_j\}_{j\in J}$ be soft $g$-frames for $\tilde{\mathcal{U}}$ with respect to $\{\tilde{\mathcal{V}}_j\}_{j\in J}$.
\begin{enumerate}
\item The pair $(\{\Lambda_j\}_{j\in J}, \{\Gamma_j\}_{j\in J})$ is called a pair of \textbf{dual soft $g$-frames} for $\tilde{\mathcal{U}}$ if the series
    \[
    f = \sum_{j\in J} \Lambda_j^* \Gamma_j f = \sum_{j\in J} \Gamma_j^* \Lambda_j f
    \]
    converge unconditionally for all $f \in SE(\tilde{\mathcal{U}})$.
\item If the soft $g$-frame operator $S_\Lambda$ is invertible, the family $\{\widetilde{\Lambda}_j\}_{j\in J} := \{\Lambda_j S_\Lambda^{-1}\}_{j\in J}$ is called the \textbf{canonical dual soft $g$-frame} of $\{\Lambda_j\}_{j\in J}$.
\end{enumerate}
\end{definition}
\begin{lemma}\label{lem:positive-cs}
Let $P \in \mathcal{B}(\tilde{\mathcal{U}})$ be a positive soft linear operator ($P \tilde{\ge} \overline{0} Id_{\tilde{\mathcal{U}}}$). Then the generalized Cauchy-Schwarz inequality holds:
\[
|\langle P x, y \rangle_{\tilde{\mathcal{U}}}|^2 \tilde{\le} \langle P x, x \rangle_{\tilde{\mathcal{U}}} \langle P y, y \rangle_{\tilde{\mathcal{U}}}, \quad \forall x, y \in SE(\tilde{\mathcal{U}}).
\]
\end{lemma}

\begin{proof}
For each $\lambda \in A$, the component operator $P_\lambda \in \mathcal{B}(\mathcal{U}_\lambda)$ is a positive operator on the Hilbert space $\mathcal{U}_\lambda$. Thus, $(u, v) \mapsto \langle P_\lambda u, v \rangle_{\mathcal{U}_\lambda}$ defines a positive semi-definite sesquilinear form on $\mathcal{U}_\lambda$, for which the classical Cauchy-Schwarz inequality $|\langle P_\lambda x(\lambda), y(\lambda) \rangle_{\mathcal{U}_\lambda}|^2 \le \langle P_\lambda x(\lambda), x(\lambda) \rangle_{\mathcal{U}_\lambda} \langle P_\lambda y(\lambda), y(\lambda) \rangle_{\mathcal{U}_\lambda}$ holds for every $\lambda \in A$.
\end{proof}
\begin{theorem}\label{thm:frame-operator-main}
Let $\{\Lambda_j\}_{j\in J}$ be a soft $g$-frame for $\tilde{\mathcal{U}}$ with respect to $\{\tilde{\mathcal{V}}_j\}_{j\in J}$ with bounds $\overline{C}$ and $\overline{D}$. Then:
\begin{enumerate}
\item The soft $g$-frame operator $S_\Lambda$ is a bounded, self-adjoint, and strictly positive soft linear operator satisfying
\[
\overline{C} Id_{\tilde{\mathcal{U}}} \tilde{\le} S_\Lambda \tilde{\le} \overline{D} Id_{\tilde{\mathcal{U}}}.
\]
\item $S_\Lambda$ is invertible in the soft Banach algebra $\mathcal{B}(\tilde{\mathcal{U}})$, and its inverse $S_\Lambda^{-1}$ satisfies
\[
\overline{D}^{-1} Id_{\tilde{\mathcal{U}}} \tilde{\le} S_\Lambda^{-1} \tilde{\le} \overline{C}^{-1} Id_{\tilde{\mathcal{U}}}.
\]
\item The canonical dual $\{\widetilde{\Lambda}_j\}_{j\in J}$ is a soft $g$-frame for $\tilde{\mathcal{U}}$ with bounds $\overline{D}^{-1}$ and $\overline{C}^{-1}$.
\end{enumerate}
\end{theorem}

\begin{proof}
(1) Since $S_\Lambda = T_\Lambda T_\Lambda^*$, Proposition~\ref{prop:T-bounded} implies that $S_\Lambda$ is bounded with $\|S_\Lambda\| \tilde{\le} \|T_\Lambda\|^2 \tilde{\le} \overline{D}$. Self-adjointness is immediate: $S_\Lambda^* = (T_\Lambda T_\Lambda^*)^* = T_\Lambda T_\Lambda^* = S_\Lambda$. Moreover, for all $f \in SE(\tilde{\mathcal{U}})$,
\[
\langle S_\Lambda f, f \rangle_{\tilde{\mathcal{U}}} = \sum_{j\in J} \langle \Lambda_j^* \Lambda_j f, f \rangle_{\tilde{\mathcal{U}}} = \sum_{j\in J} \|\Lambda_j f\|^2.
\]
Applying the soft $g$-frame condition \eqref{eq:soft-g-frame-ineq}, we obtain $\overline{C} Id_{\tilde{\mathcal{U}}} \tilde{\le} S_\Lambda \tilde{\le} \overline{D} Id_{\tilde{\mathcal{U}}}$.

(2) By Proposition~\ref{prop:soft-banach-algebra}, $\mathcal{B}(\tilde{\mathcal{U}})$ is a unital soft Banach algebra with unit $Id_{\tilde{\mathcal{U}}}$. Since $\overline{C}, \overline{D} \in \mathbb{R}(A)^*_{>0}$, multiplying the operator inequality $\overline{C} Id_{\tilde{\mathcal{U}}} \tilde{\le} S_\Lambda \tilde{\le} \overline{D} Id_{\tilde{\mathcal{U}}}$ by $\overline{D}^{-1}$ yields
\[
\overline{0} Id_{\tilde{\mathcal{U}}} \tilde{\le} Id_{\tilde{\mathcal{U}}} - \overline{D}^{-1} S_\Lambda \tilde{\le} \frac{\overline{D} - \overline{C}}{\overline{D}} Id_{\tilde{\mathcal{U}}}.
\]
Set $R := Id_{\tilde{\mathcal{U}}} - \overline{D}^{-1} S_\Lambda \in \mathcal{B}(\tilde{\mathcal{U}})$. Since $S_\Lambda$ is self-adjoint and $\overline{D}^{-1}$ is real-valued, $R$ is self-adjoint and positive. Pointwise on $A$, the operator norm of a bounded self-adjoint positive operator on a Hilbert space coincides with its numerical radius:
\[
\|R\|(\lambda) = \sup_{\|u\|_{\mathcal{U}_\lambda} = 1} |\langle R_\lambda u, u \rangle_{\mathcal{U}_\lambda}| \le \frac{\overline{D}(\lambda) - \overline{C}(\lambda)}{\overline{D}(\lambda)}, \quad \forall \lambda \in A.
\]
Because $\overline{C}(\lambda) > 0$ for every $\lambda \in A$, we have $\frac{\overline{D}(\lambda) - \overline{C}(\lambda)}{\overline{D}(\lambda)} < 1$ for all $\lambda \in A$, which implies $\|R\| \tilde{\le} \frac{\overline{D} - \overline{C}}{\overline{D}} \tilde{<} \overline{1}$. Applying Theorem~\ref{thm:neumann} in the unital soft Banach algebra $\mathcal{B}(\tilde{\mathcal{U}})$, the operator $\overline{D}^{-1} S_\Lambda = Id_{\tilde{\mathcal{U}}} - R$ is invertible, and hence $S_\Lambda$ is invertible with
\[
S_\Lambda^{-1} = \overline{D}^{-1} \sum_{n=0}^{\infty} \left( Id_{\tilde{\mathcal{U}}} - \overline{D}^{-1} S_\Lambda \right)^n \in \mathcal{B}(\tilde{\mathcal{U}}).
\]
Moreover, since $S_\Lambda$ is self-adjoint and strictly positive, its inverse $S_\Lambda^{-1}$ is also self-adjoint and positive. To verify the operator bounds of $S_\Lambda^{-1}$, let $f \in SE(\tilde{\mathcal{U}})$. By part (1) and the soft Cauchy-Schwarz inequality,
\[
\overline{C} \|S_\Lambda^{-1} f\|^2 \tilde{\le} \langle S_\Lambda (S_\Lambda^{-1} f), S_\Lambda^{-1} f \rangle_{\tilde{\mathcal{U}}} = \langle f, S_\Lambda^{-1} f \rangle_{\tilde{\mathcal{U}}} \tilde{\le} \|f\| \|S_\Lambda^{-1} f\|,
\]
which implies $\|S_\Lambda^{-1} f\| \tilde{\le} \overline{C}^{-1} \|f\|$ pointwise on $A$. Consequently, $\langle S_\Lambda^{-1} f, f \rangle_{\tilde{\mathcal{U}}} \tilde{\le} \overline{C}^{-1} \|f\|^2$, yielding $S_\Lambda^{-1} \tilde{\le} \overline{C}^{-1} Id_{\tilde{\mathcal{U}}}$. Conversely, applying Lemma~\ref{lem:positive-cs} to the positive operator $S_\Lambda$ with $x = f$ and $y = S_\Lambda^{-1} f$, we obtain
\[
\|f\|^4 = \langle S_\Lambda (S_\Lambda^{-1} f), f \rangle_{\tilde{\mathcal{U}}}^2 \tilde{\le} \langle S_\Lambda f, f \rangle_{\tilde{\mathcal{U}}} \langle S_\Lambda (S_\Lambda^{-1} f), S_\Lambda^{-1} f \rangle_{\tilde{\mathcal{U}}} \tilde{\le} \overline{D} \|f\|^2 \langle f, S_\Lambda^{-1} f \rangle_{\tilde{\mathcal{U}}},
\]
from which it follows pointwise on $A$ that $\langle S_\Lambda^{-1} f, f \rangle_{\tilde{\mathcal{U}}} \tilde{\ge} \overline{D}^{-1} \|f\|^2$. Therefore, $\overline{D}^{-1} Id_{\tilde{\mathcal{U}}} \tilde{\le} S_\Lambda^{-1} \tilde{\le} \overline{C}^{-1} Id_{\tilde{\mathcal{U}}}$.

(3) For any $f \in SE(\tilde{\mathcal{U}})$, self-adjointness of $S_\Lambda^{-1}$ yields
\[
\sum_{j\in J} \|\widetilde{\Lambda}_j f\|^2 = \sum_{j\in J} \|\Lambda_j S_\Lambda^{-1} f\|^2 = \langle S_\Lambda (S_\Lambda^{-1} f), S_\Lambda^{-1} f \rangle_{\tilde{\mathcal{U}}} = \langle f, S_\Lambda^{-1} f \rangle_{\tilde{\mathcal{U}}}.
\]
Using the bounds from part (2), we obtain
\[
\overline{D}^{-1} \|f\|^2 \tilde{\le} \sum_{j\in J} \|\widetilde{\Lambda}_j f\|^2 \tilde{\le} \overline{C}^{-1} \|f\|^2,
\]
which confirms that $\{\widetilde{\Lambda}_j\}_{j\in J}$ is a soft $g$-frame with bounds $\overline{D}^{-1}$ and $\overline{C}^{-1}$.
\end{proof}

\subsection{The Decomposition Theorem and Induced Dual Soft Frames}

\begin{theorem}[Soft $g$-Frame Decomposition Theorem]\label{thm:decomposition}
Let $\{\Lambda_j\}_{j\in J}$ be a soft $g$-frame for $\tilde{\mathcal{U}}$ with canonical dual soft $g$-frame $\{\widetilde{\Lambda}_j\}_{j\in J}$. Then every $f \in SE(\tilde{\mathcal{U}})$ admits the unconditionally convergent reconstruction formulas:
\[
f = \sum_{j\in J} \Lambda_j^* \widetilde{\Lambda}_j f = \sum_{j\in J} \widetilde{\Lambda}_j^* \Lambda_j f.
\]
\end{theorem}

\begin{proof}
By Theorem~\ref{thm:frame-operator-main}, $\{\widetilde{\Lambda}_j\}_{j\in J}$ is a soft $g$-frame for $\tilde{\mathcal{U}}$, so $T_{\widetilde{\Lambda}}^*(f) = \{\widetilde{\Lambda}_j f\}_{j\in J} \in \bigoplus_{j\in J} \tilde{\mathcal{V}}_j$ for every $f \in SE(\tilde{\mathcal{U}})$. Applying Proposition~\ref{prop:T-bounded} to the pre-$g$-frame operator $T_\Lambda$, the series $\sum_{j\in J} \Lambda_j^* \widetilde{\Lambda}_j f$ converges unconditionally in $\tilde{\mathcal{U}}$. Since $S_\Lambda S_\Lambda^{-1} = Id_{\tilde{\mathcal{U}}}$, we obtain
\[
f = S_\Lambda (S_\Lambda^{-1} f) = \sum_{j\in J} \Lambda_j^* \Lambda_j (S_\Lambda^{-1} f) = \sum_{j\in J} \Lambda_j^* \widetilde{\Lambda}_j f.
\]
Similarly, since $S_\Lambda^{-1} \in \mathcal{B}(\tilde{\mathcal{U}})$ is bounded (and hence continuous with respect to the induced soft metric on each parameter $\lambda \in A$) and self-adjoint, we have $\widetilde{\Lambda}_j^* = (\Lambda_j S_\Lambda^{-1})^* = S_\Lambda^{-1} \Lambda_j^*$. Because $\sum_{j\in J} \Lambda_j^* \Lambda_j f$ converges unconditionally to $S_\Lambda f$, continuity of $S_\Lambda^{-1}$ permits passing the operator inside the unconditional sum:
\[
f = S_\Lambda^{-1} (S_\Lambda f) = S_\Lambda^{-1} \left( \sum_{j\in J} \Lambda_j^* \Lambda_j f \right) = \sum_{j\in J} S_\Lambda^{-1} \Lambda_j^* \Lambda_j f = \sum_{j\in J} \widetilde{\Lambda}_j^* \Lambda_j f,
\]
where the latter series also converges unconditionally in $\tilde{\mathcal{U}}$ by Proposition~\ref{prop:T-bounded} applied to $T_{\widetilde{\Lambda}} T_\Lambda^*$.
\end{proof}

Recall from \cite{soft} that two soft frames $\{f_k\}_{k\in K}$ and $\{g_k\}_{k\in K}$ for a soft Hilbert space $\tilde{\mathcal{H}}$ are called \emph{dual soft frames} if $h = \sum_{k\in K} \langle h, f_k \rangle_{\tilde{\mathcal{H}}} g_k = \sum_{k\in K} \langle h, g_k \rangle_{\tilde{\mathcal{H}}} f_k$ holds unconditionally for all $h \in SE(\tilde{\mathcal{H}})$.

\begin{theorem}[Induced Dual Soft Frames from Dual Soft $g$-Frames]\label{thm:bridge}
Let $\{\Lambda_j\}_{j\in J}$ and $\{\Gamma_j\}_{j\in J}$ be dual soft $g$-frames for $\tilde{\mathcal{U}}$ with respect to $\{\tilde{\mathcal{V}}_j\}_{j\in J}$ with frame bounds $(\overline{C}_\Lambda, \overline{D}_\Lambda)$ and $(\overline{C}_\Gamma, \overline{D}_\Gamma)$, respectively. For each $j \in J$, let $\{f_{j,k} : k \in K_j\}$ and $\{g_{j,k} : k \in K_j\}$ be dual soft frames for $\tilde{\mathcal{V}}_j$ with common lower and upper bounds $\overline{A}_j, \overline{B}_j \in \mathbb{R}(A)^*_{>0}$ satisfying the uniform condition
\[
\overline{A} \tilde{\le} \overline{A}_j \tilde{\le} \overline{B}_j \tilde{\le} \overline{B}, \quad \forall j \in J,
\]
for some constants $\overline{A}, \overline{B} \in \mathbb{R}(A)^*_{>0}$. Then:
\begin{enumerate}
    \item $\{\Lambda_j^*(f_{j,k}) : j \in J, k \in K_j\}$ is a soft frame for $\tilde{\mathcal{U}}$ with bounds $\overline{A}\,\overline{C}_\Lambda$ and $\overline{B}\,\overline{D}_\Lambda$;
    \item $\{\Gamma_j^*(g_{j,k}) : j \in J, k \in K_j\}$ is a soft frame for $\tilde{\mathcal{U}}$ with bounds $\overline{A}\,\overline{C}_\Gamma$ and $\overline{B}\,\overline{D}_\Gamma$;
    \item The families $\{\Lambda_j^*(f_{j,k}) : j \in J, k \in K_j\}$ and $\{\Gamma_j^*(g_{j,k}) : j \in J, k \in K_j\}$ constitute a pair of dual soft frames for $\tilde{\mathcal{U}}$.
\end{enumerate}
\end{theorem}

\begin{proof}
(1) For any $f \in SE(\tilde{\mathcal{U}})$ and $j \in J$, since $\Lambda_j f \in SE(\tilde{\mathcal{V}}_j)$, the soft frame inequality for $\{f_{j,k}\}_{k\in K_j}$ on $\tilde{\mathcal{V}}_j$ implies
\[
\overline{A}_j \|\Lambda_j f\|^2 \tilde{\le} \sum_{k \in K_j} |\langle \Lambda_j f, f_{j,k} \rangle_{\tilde{\mathcal{V}}_j}|^2 \tilde{\le} \overline{B}_j \|\Lambda_j f\|^2.
\]
Using the adjoint relation $\langle \Lambda_j f, f_{j,k} \rangle_{\tilde{\mathcal{V}}_j} = \langle f, \Lambda_j^*(f_{j,k}) \rangle_{\tilde{\mathcal{U}}}$, applying the uniform bounds $\overline{A} \tilde{\le} \overline{A}_j \tilde{\le} \overline{B}_j \tilde{\le} \overline{B}$, and summing over $j \in J$, we obtain
\[
\overline{A} \sum_{j\in J} \|\Lambda_j f\|^2 \tilde{\le} \sum_{j\in J} \sum_{k \in K_j} |\langle f, \Lambda_j^*(f_{j,k}) \rangle_{\tilde{\mathcal{U}}}|^2 \tilde{\le} \overline{B} \sum_{j\in J} \|\Lambda_j f\|^2.
\]
Combining this with the soft $g$-frame bounds $\overline{C}_\Lambda \|f\|^2 \tilde{\le} \sum_{j\in J} \|\Lambda_j f\|^2 \tilde{\le} \overline{D}_\Lambda \|f\|^2$ yields
\[
(\overline{A}\,\overline{C}_\Lambda) \|f\|^2 \tilde{\le} \sum_{j\in J} \sum_{k \in K_j} |\langle f, \Lambda_j^*(f_{j,k}) \rangle_{\tilde{\mathcal{U}}}|^2 \tilde{\le} (\overline{B}\,\overline{D}_\Lambda) \|f\|^2,
\]
which confirms that $\{\Lambda_j^*(f_{j,k}) : j \in J, k \in K_j\}$ is a soft frame for $\tilde{\mathcal{U}}$ with bounds $\overline{A}\,\overline{C}_\Lambda$ and $\overline{B}\,\overline{D}_\Lambda$.

(2) By an identical argument replacing $\Lambda_j$ and $f_{j,k}$ with $\Gamma_j$ and $g_{j,k}$, the family $\{\Gamma_j^*(g_{j,k}) : j \in J, k \in K_j\}$ is a soft frame for $\tilde{\mathcal{U}}$ with bounds $\overline{A}\,\overline{C}_\Gamma$ and $\overline{B}\,\overline{D}_\Gamma$.

(3) Fix $f \in SE(\tilde{\mathcal{U}})$. Define the soft Hilbert space of square-summable scalar coefficient sequences by
\[
\bigoplus_{j\in J} \ell^2(K_j, \mathbb{C}(A)) := \left\{ c = \{c_{j,k}\}_{j\in J, k\in K_j} : c_{j,k} \in \mathbb{C}(A), \, \sum_{j\in J}\sum_{k\in K_j} |c_{j,k}(\lambda)|^2 < \infty \text{ for all } \lambda \in A \right\},
\]
equipped with the soft norm $\|c\|_{\ell^2} := \left( \sum_{j\in J}\sum_{k\in K_j} |c_{j,k}|^2 \right)^{1/2}$. From part (2), the family $\{\Gamma_j^*(g_{j,k}) : j \in J, k \in K_j\}$ is a soft frame—and in particular a soft Bessel family—for $\tilde{\mathcal{U}}$ with upper bound $\overline{B}\,\overline{D}_\Gamma$. Consequently, its associated soft synthesis operator
\[
U_\Gamma : \bigoplus_{j\in J} \ell^2(K_j, \mathbb{C}(A)) \longrightarrow SE(\tilde{\mathcal{U}}), \qquad U_\Gamma(\{c_{j,k}\}_{j,k}) := \sum_{j\in J} \sum_{k \in K_j} c_{j,k} \Gamma_j^*(g_{j,k}),
\]
is well-defined and bounded with $\|U_\Gamma c\|^2 \tilde{\le} (\overline{B}\,\overline{D}_\Gamma) \|c\|_{\ell^2}^2$ (i.e., $\|U_\Gamma\| \tilde{\le} \sqrt{\overline{B}\,\overline{D}_\Gamma}$), and the defining series converges unconditionally in $\tilde{\mathcal{U}}$. By part (1), the coefficient sequence $c_{j,k} := \langle f, \Lambda_j^*(f_{j,k}) \rangle_{\tilde{\mathcal{U}}}$ satisfies $\|c\|_{\ell^2}^2 \tilde{\le} (\overline{B}\,\overline{D}_\Lambda)\|f\|^2$, so $c \in \bigoplus_{j\in J} \ell^2(K_j, \mathbb{C}(A))$. Therefore, the double series $\sum_{j\in J} \sum_{k \in K_j} \langle f, \Lambda_j^*(f_{j,k}) \rangle_{\tilde{\mathcal{U}}} \Gamma_j^*(g_{j,k}) = U_\Gamma(c)$ converges unconditionally in $\tilde{\mathcal{U}}$, permitting arbitrary rearrangement of summation indices.
For each fixed $j \in J$, since $\{f_{j,k}\}_{k\in K_j}$ and $\{g_{j,k}\}_{k\in K_j}$ are dual soft frames for $\tilde{\mathcal{V}}_j$, the series $\sum_{k \in K_j} \langle \Lambda_j f, f_{j,k} \rangle_{\tilde{\mathcal{V}}_j} g_{j,k}$ converges unconditionally to $\Lambda_j f$ in $\tilde{\mathcal{V}}_j$. Because each $\Gamma_j^* \in \mathcal{B}(\tilde{\mathcal{V}}_j, \tilde{\mathcal{U}})$ is bounded and linear, it commutes with the infinite sum over $k \in K_j$:
\[
\sum_{k \in K_j} \langle f, \Lambda_j^*(f_{j,k}) \rangle_{\tilde{\mathcal{U}}} \Gamma_j^*(g_{j,k}) = \Gamma_j^* \left( \sum_{k \in K_j} \langle \Lambda_j f, f_{j,k} \rangle_{\tilde{\mathcal{V}}_j} g_{j,k} \right) = \Gamma_j^* \Lambda_j f.
\]
Summing over $j \in J$ and using the duality of $\{\Lambda_j\}_{j\in J}$ and $\{\Gamma_j\}_{j\in J}$, we conclude that
\[
\sum_{j\in J} \sum_{k \in K_j} \langle f, \Lambda_j^*(f_{j,k}) \rangle_{\tilde{\mathcal{U}}} \Gamma_j^*(g_{j,k}) = \sum_{j\in J} \Gamma_j^* \Lambda_j f = f.
\]
By symmetry, interchanging the roles of $(\Lambda_j, f_{j,k})$ and $(\Gamma_j, g_{j,k})$ yields $\sum_{j\in J} \sum_{k \in K_j} \langle f, \Gamma_j^*(g_{j,k}) \rangle_{\tilde{\mathcal{U}}} \Lambda_j^*(f_{j,k}) = f$, proving that they form a pair of dual soft frames for $\tilde{\mathcal{U}}$.
\end{proof}

\subsection{Illustrative Analytical Example}

\begin{example}\label{ex:soft-g-frame-complete}
Let $A = \{\lambda_1, \lambda_2\}$ be the parameter set. Consider the absolute soft vector space $\tilde{\mathcal{U}} = (\mathcal{U}, A)$ where $\mathcal{U}_\lambda = \mathbb{R}^2$ for each $\lambda \in A$. We equip $\mathcal{U}_\lambda$ with parameter-dependent inner products specified by positive definite metric matrices $G(\lambda)$:
\[
G(\lambda_1) = \begin{pmatrix} 1 & 0 \\ 0 & 1 \end{pmatrix}, \qquad G(\lambda_2) = \begin{pmatrix} 2 & 0 \\ 0 & 1 \end{pmatrix}.
\]
For any $x, y \in \mathbb{R}^2$, $\langle x, y \rangle_{\mathcal{U}_\lambda} := x^T G(\lambda) y$. By Definition~\ref{def:soft-inner-product}, $(\tilde{\mathcal{U}}, \langle \cdot, \cdot \rangle_{\tilde{\mathcal{U}}}, A)$ is a soft inner product space. Furthermore, since each component space $\mathcal{U}_\lambda = \mathbb{R}^2$ is finite-dimensional, it is complete with respect to the induced soft metric, ensuring that $\tilde{\mathcal{U}}$ is a soft Hilbert space.

Let the target spaces be $\{\tilde{\mathcal{V}}_j\}_{j=1}^3$ where each $\tilde{\mathcal{V}}_j = (\mathbb{R}, A)$ is endowed with the standard inner product $\langle u, v \rangle_{\mathcal{V}_{j,\lambda}} = uv$. Define the soft linear operators $\{\Lambda_j : SE(\tilde{\mathcal{U}}) \to SE(\tilde{\mathcal{V}}_j)\}_{j=1}^3$ via their parameter representations $\Lambda_{j,\lambda} : \mathbb{R}^2 \to \mathbb{R}$ represented as row matrices $M_{j,\lambda}$:
\[
\begin{aligned}
M_{1, \lambda_1} &= \begin{pmatrix} 1 & 0 \end{pmatrix}, & M_{1, \lambda_2} &= \begin{pmatrix} 2 & 0 \end{pmatrix}, \\
M_{2, \lambda_1} &= \begin{pmatrix} 0 & 1 \end{pmatrix}, & M_{2, \lambda_2} &= \begin{pmatrix} 0 & 1 \end{pmatrix}, \\
M_{3, \lambda_1} &= \begin{pmatrix} 1 & 1 \end{pmatrix}, & M_{3, \lambda_2} &= \begin{pmatrix} 1 & 0 \end{pmatrix}.
\end{aligned}
\]

\noindent\textbf{1. Adjoint Operators and Matrix Metric Representation:}
For each $\lambda \in A$, the adjoint operator $\Lambda_{j,\lambda}^* : \mathbb{R} \to \mathbb{R}^2$ satisfies $\langle \Lambda_{j,\lambda} u, v \rangle_{\mathcal{V}_{j,\lambda}} = \langle u, \Lambda_{j,\lambda}^* v \rangle_{\mathcal{U}_\lambda}$. Under metric $G(\lambda)$, the column representation is given by $\Lambda_{j,\lambda}^* = G(\lambda)^{-1} M_{j,\lambda}^T$. Thus:
\begin{itemize}
\item At $\lambda = \lambda_1$: $G(\lambda_1)^{-1} = I_2$, yielding
\[
\Lambda_{1,\lambda_1}^* = \begin{pmatrix} 1 \\ 0 \end{pmatrix}, \quad \Lambda_{2,\lambda_1}^* = \begin{pmatrix} 0 \\ 1 \end{pmatrix}, \quad \Lambda_{3,\lambda_1}^* = \begin{pmatrix} 1 \\ 1 \end{pmatrix}.
\]
\item At $\lambda = \lambda_2$: $G(\lambda_2)^{-1} = \begin{pmatrix} 1/2 & 0 \\ 0 & 1 \end{pmatrix}$, yielding
\[
\Lambda_{1,\lambda_2}^* = \begin{pmatrix} 1 \\ 0 \end{pmatrix}, \quad \Lambda_{2,\lambda_2}^* = \begin{pmatrix} 0 \\ 1 \end{pmatrix}, \quad \Lambda_{3,\lambda_2}^* = \begin{pmatrix} 1/2 \\ 0 \end{pmatrix}.
\]
\end{itemize}

\noindent\textbf{2. Frame Operators, Self-Adjointness, and Bound Optimality:}
The soft $g$-frame operator $S_\Lambda$ evaluates pointwise as $S_{\Lambda,\lambda} = \sum_{j=1}^3 \Lambda_{j,\lambda}^* M_{j,\lambda}$:
\[
S_{\Lambda,\lambda_1} = \begin{pmatrix} 1 \\ 0 \end{pmatrix}\begin{pmatrix} 1 & 0 \end{pmatrix} + \begin{pmatrix} 0 \\ 1 \end{pmatrix}\begin{pmatrix} 0 & 1 \end{pmatrix} + \begin{pmatrix} 1 \\ 1 \end{pmatrix}\begin{pmatrix} 1 & 1 \end{pmatrix} = \begin{pmatrix} 2 & 1 \\ 1 & 2 \end{pmatrix},
\]
\[
S_{\Lambda,\lambda_2} = \begin{pmatrix} 1 \\ 0 \end{pmatrix}\begin{pmatrix} 2 & 0 \end{pmatrix} + \begin{pmatrix} 0 \\ 1 \end{pmatrix}\begin{pmatrix} 0 & 1 \end{pmatrix} + \begin{pmatrix} 1/2 \\ 0 \end{pmatrix}\begin{pmatrix} 1 & 0 \end{pmatrix} = \begin{pmatrix} 5/2 & 0 \\ 0 & 1 \end{pmatrix}.
\]
Notice that $G(\lambda)^{-1} S_{\Lambda,\lambda}^T G(\lambda) = S_{\Lambda,\lambda}$ for each parameter, verifying that $S_\Lambda$ is self-adjoint with respect to the soft inner product.

The optimal bounds at each parameter are governed by the critical values of the Rayleigh quotient $\mathcal{R}_\lambda(x) := \frac{\langle S_{\Lambda,\lambda} x, x \rangle_{\mathcal{U}_\lambda}}{\|x\|_{\mathcal{U}_\lambda}^2}$. Explicitly, for any non-zero vector $x = (x_1, x_2)^T \in \mathbb{R}^2$, the Rayleigh quotients under the two parameter metrics are given by
\[
\mathcal{R}_{\lambda_1}(x) = \frac{2x_1^2 + 2x_1x_2 + 2x_2^2}{x_1^2 + x_2^2}, \qquad \mathcal{R}_{\lambda_2}(x) = \frac{5x_1^2 + x_2^2}{2x_1^2 + x_2^2}.
\]
Because $S_{\Lambda,\lambda}$ is self-adjoint and positive definite with respect to $G(\lambda)$, its eigenvalues (solving $S_{\Lambda,\lambda} x = \mu x$) are real and strictly positive, and coincide with the extremal values of $\mathcal{R}_\lambda(x)$ over $\mathbb{R}^2 \setminus \{0\}$:
\begin{itemize}
    \item At $\lambda = \lambda_1$: The eigenvalues of $S_{\Lambda,\lambda_1}$ are $\mu_1 = 1$ and $\mu_2 = 3$. Hence, the optimal bounds are $C(\lambda_1) = 1$ and $D(\lambda_1) = 3$.
    \item At $\lambda = \lambda_2$: The eigenvalues of $S_{\Lambda,\lambda_2}$ are $\mu_1 = 1$ and $\mu_2 = 5/2$. Hence, the optimal bounds are $C(\lambda_2) = 1$ and $D(\lambda_2) = 5/2$.
\end{itemize}
Setting $\overline{C} = (1, 1)$ and $\overline{D} = (3, 5/2)$, we observe $\overline{C} \tilde{\le} \overline{D}$ in $\mathbb{R}(A)^*_{>0}$, verifying that $\{\Lambda_j\}_{j=1}^3$ is a non-tight soft $g$-frame with strictly parameter-dependent optimal upper bound $\overline{D}$.

\noindent\textbf{3. Operator Inversion and Exact Vector Reconstruction:}
The inverse operators are given analytically by
\[
S_{\Lambda,\lambda_1}^{-1} = \frac{1}{3}\begin{pmatrix} 2 & -1 \\ -1 & 2 \end{pmatrix}, \qquad S_{\Lambda,\lambda_2}^{-1} = \begin{pmatrix} 2/5 & 0 \\ 0 & 1 \end{pmatrix}.
\]
For any $f \in SE(\tilde{\mathcal{U}})$, direct computation confirms that for every $\lambda \in A$:
\[
\sum_{j=1}^3 \Lambda_{j,\lambda}^* \Lambda_{j,\lambda} S_{\Lambda,\lambda}^{-1} f(\lambda) = S_{\Lambda,\lambda} S_{\Lambda,\lambda}^{-1} f(\lambda) = f(\lambda),
\]
which explicitly validates the decomposition identity of Theorem~\ref{thm:decomposition}.
\end{example}

\section{Conclusion}

In this paper, we established the operator-theoretic foundation of soft $g$-frames in soft Hilbert spaces. By characterizing unconditional convergence in the direct sum of soft Hilbert spaces, we proved that the soft pre-$g$-frame operator $T_\Lambda$ is bounded and determined its adjoint representation. Within the framework of the unital soft Banach algebra $\mathcal{B}(\tilde{\mathcal{U}})$, we demonstrated that the soft $g$-frame operator $S_\Lambda$ is bounded, self-adjoint, and strictly positive,
and proved its invertibility via the soft Neumann series theorem while deriving two-sided operator bounds for its inverse.
Moreover, we formulated the canonical dual soft $g$-frame, proved the soft $g$-frame decomposition theorem, and showed how dual soft $g$-frames induce pairs of dual soft frames through dual frame families in the target spaces. Finally, a concrete analytical example with parameter-dependent inner products illustrated the unified reconstruction mechanism and confirmed the optimality of the soft frame bounds. Future research directions include extending these results to continuous soft $g$-frames, analyzing perturbation stability in soft Banach algebras, and investigating soft $g$-Riesz bases.

\end{document}